\documentclass[runningheads]{llncs}
\usepackage[T1]{fontenc}
\usepackage{graphicx}
\usepackage{hyperref}
\usepackage{color}

\usepackage{amsmath, amssymb}
\usepackage{booktabs}
\usepackage{todonotes}
\begin{document}
%
%\title{Improving upon LLM-driven discovery: \\                Global Optimization for \\ Hard Packing Problems}
 \title{Global Optimization for \\ Combinatorial Geometry Problems \\ Revisited in the Era of LLMs}
% \title{Global Optimization for Combinatorial Geometry Problems in the Era of LLM-based Code Generators}
% \title{A Global Optimization Perspective on Hard Packing Problems                 Highlighted by Recent LLM-Based Innovation}
% \title{How Global Optimization Can Solve Hard Packing Problems                 Recently Highlighted by LLM Approaches}

\titlerunning{Global Optimization for Combinatorial Geometry Problems}
\author{Timo Berthold\inst{1,2} \and
Dominik Kamp\inst{3} \and
Gioni Mexi\inst{3} \and
Sebastian Pokutta\inst{2,3} \and
Imre P\'olik\inst{4}}
\authorrunning{T. Berthold et al.}
\institute{Fair Isaac Deutschland GmbH, Germany, \email{timoberthold@fico.com} \and
Technische Universit\"at Berlin, Institut f\"ur Mathematik, Germany \and
Zuse Institute Berlin, Germany \and
Fair Isaac Corporation, USA, \email{imre@polik.net}}
\maketitle              % typeset the header of the contribution
\begin{abstract}
Recent progress in LLM-driven algorithm discovery, exemplified by DeepMind's AlphaEvolve, has produced new best-known solutions for a range of hard geometric and combinatorial problems. This raises a natural question: to what extent can modern off-the-shelf global optimization solvers match such results when the problems are formulated directly as nonlinear optimization problems (NLPs)?

We revisit a subset of problems from the AlphaEvolve benchmark suite and evaluate straightforward NLP formulations with two state-of-the-art solvers, the commercial FICO Xpress and the open-source SCIP. Without any solver modifications, both solvers reproduce, and in several cases improve upon, the best solutions previously reported in the literature, including the recent LLM-driven discoveries. Our results not only highlight the maturity of generic NLP technology and its ability to tackle nonlinear mathematical problems that were out of reach for general-purpose solvers only a decade ago, but also position global NLP solvers as powerful tools that may be exploited within LLM-driven algorithm discovery.
\end{abstract}
\section{Introduction}

The rapid progress in global optimization technology over the past decade has substantially expanded the range of nonlinear, nonconvex problems that can be solved to proven global optimality, or at least to very high-quality solutions with reliable dual bounds.
State-of-the-art academic solvers like SCIP~\cite{SCIPOptSuite10} and commercial solvers like FICO\textsuperscript{®} Xpress~\cite{BelottiBertholdGallyGottwaldPolik2025} combine spatial branch-and-bound, automatic linearization and convexification, sophisticated presolving, and increasingly powerful primal heuristics~\cite{Berthold2014,berthold2018computational,Berthold_Lodi_Salvagnin_2025}, enabling them to solve instances that would have been considered computationally prohibitive only a few years ago.

At the same time, recent developments in algorithm design based on Large Language Models (LLMs) have drawn renewed attention to long-standing geometric and combinatorial problems that can be formulated as nonlinear optimization models.
Specifically, DeepMind presented the AlphaEvolve framework~\cite{georgiev2025mathematical,novikov2025alphaevolve}, which uses LLM-generated code in an evolutionary search to produce high-quality solutions for an extensive set of mathematical problems, including variants of circle packing, hexagon packing, and minimum-distance configurations of points.
These breakthrough results raise a natural question: to what extent can state-of-the-art global optimization solvers match, or even surpass, such automatically discovered algorithms on challenging nonlinear problems?

More broadly, the last few years have seen rapid progress in LLM-driven discovery workflows that couple generative models with structured search and automated evaluation or verification.
FunSearch combines an LLM with evolutionary program search and task-specific evaluators, enabling improvements for several discrete mathematical and algorithmic problems~\cite{romeraparedes2024funsearch}.
Complementary advances include AlphaDev, which used learning-based search to rediscover and improve low-level algorithms such as sorting routines~\cite{Mankowitz2023}, and AlphaGeometry, which combines neural generation with symbolic reasoning to solve olympiad-level geometry problems without human demonstrations~\cite{trinh2024solving}.
Related LLM+evolution approaches have also been used to design effective heuristics for combinatorial optimization~\cite{liu2024example}.
Taken together, these advances motivate revisiting classical mathematical optimization as a competitive and reliable baseline on the same challenging benchmarks.

In this work we revisit three problems from the AlphaEvolve benchmark suite and study them through the lens of Nonlinear Programming (NLP).
An \emph{NLP} is an optimization problem minimizing a nonlinear objective function over a feasible set defined by nonlinear constraints on continuous variables:
\begin{align}
    \min\{f(x)\mid g_k(x)\leq 0,\forall k\in\mathcal{K}, \ell\leq x\leq u\},
    \label{nlp_def}
\end{align} where the objective $f(x)$ and all constraint functions $g_k:\mathbb{R}^{n}\rightarrow\mathbb{R}$ are factorable and all variable bounds $\ell,u\in\bar{\mathbb{R}} := \mathbb{R}\cup\{\pm \infty\}$.
The set $\mathcal{K}:=\{1,\dots,m\}$ indexes the constraints.
The circle packing problem, the minimum-distance ratio problem and the hexagon packing problem admit compact NLP models.
This makes them good showcases for the power of modern global optimization tools: these combinatorial problems can be very intuitively modelled and effectively solved in their most natural form, as nonlinear optimization instances.
We show that combined with off-the-shelf global optimization technology, these straightforward NLP formulations not only reproduce the best solutions reported in \cite{georgiev2025mathematical,novikov2025alphaevolve}, but in several cases produce significantly better solutions.
Our goal in doing so is to illustrate the power of modern general-purpose nonlinear optimization solvers rather than to perform a head-to-head comparison of solvers.

%\subsection*{Contribution}
%This paper revisits a subset of AlphaEvolve benchmarks through the lens of NLP and demonstrates that generic global optimization has matured into a strong baseline for such problems. Our main contributions are:
Next to presenting simple formulations that work and identifying some key modeling decisions that drive performance, our main contributions are:
\paragraph{Stronger solutions with unmodified solvers.}
With off-the-shelf Xpress and SCIP we obtain solutions that match, and in multiple cases improve, previously reported best-known results. Some problems were solved with only one solver, some used both. In this paper we do not make a distinction of which solver found which solution. All the solutions that our solvers produced were verified using the validation code in the AlphaEvolve repository.
\paragraph{Lessons learned.}
We reflect on the respective strengths and limitations of LLM-based approaches and nonlinear optimization, discuss modeling insights, and how global optimization software has become a complementary, industry-ready tool for challenging combinatorial problems.

These insights are discussed in Section~\ref{thats_it_folks}, while Sections~\ref{minmaxr}--\ref{hexhex} present one model each, together with the accompanying computational results. Due to lack of space we will not comment on the dual bounds and how they could be improved (such as stronger formulations or symmetry breaking), even though it is a key differentiator between ad-hoc heuristics and general optimization solvers.

\section{Minimizing the Ratio of Maximum to Minimum Distance}
\label{minmaxr}

The problem of minimizing the ratio between the maximum and minimum pairwise distances in a finite point set, called the \emph{min-max ratio problem}, has its origins in classical extremal geometry.
Bateman and Erd\H{o}s~\cite{bateman1951geometrical} formulated the question in terms of determining, for a given number of points, the configuration of points in the plane with mutual distances at least one that minimizes the diameter of the set.
This is equivalent to minimizing the ratio of largest to smallest distance when the minimum is normalized to one.
Bateman and Erd\H{o}s provided optimal solutions for $n \leq 7$ in two dimensions.
David Cantrell, see, e.g.,~\cite{friedmanPacking}, has contributed many best-known solutions in both 2D and 3D, while Audet et al.\ formulated the problem explicitly as a nonlinear optimization model and found some improving configurations for $n \leq 30$, using derivative-free optimization techniques~\cite{audet2010note}. 

\subsection{Optimization model}

In a slightly more general form, for given positive integers $n$ and $d$, the goal of the min-max ratio problem is to find $n$ distinct points in
$d$-dimensional Euclidean space $\mathbb{R}^d$ such that the ratio between the
maximum and the minimum pairwise distance is minimized. This problem seeks a configuration of points that is as evenly spaced as possible,
minimizing clustering and maximizing spatial efficiency.
Let $\mathcal{N} = \{1, 2, \ldots, n\}$ denote the set of points and 
$\mathcal{D} = \{1, 2, \ldots, d\}$ denote the set of dimensions. We define the
following decision variables:
\begin{description}
\item[$x_{i,k} \in \mathbb{R}$:] coordinate of point $i \in \mathcal{N}$ in dimension $k \in \mathcal{D}$
\item[$t_{\min} \in \mathbb{R}_+$:] minimum squared distance between any two points
\item[$t_{\max} \in \mathbb{R}_+$:] maximum squared distance between any two points
\item[$r \in \mathbb{R}_+$:] squared ratio variable, defined as $r = \frac{t_{\max}}{t_{\min}}$
\end{description}
The first observation is that it is more convenient and efficient to work with the ratio of the squared distances and avoid the square roots. For any pair of points $i, j\in\mathcal{N}$, the squared Euclidean distance is given by:
\[
d_{ij}^2 = \sum_{k \in \mathcal{D}} (x_{i,k} - x_{j,k})^2
\]
A key insight is that this problem exhibits \emph{scale-invariance}: the optimal
objective value remains unchanged if all lengths are multiplied by the same constant
factor. We can therefore normalize the problem by fixing either $t_{\min}$ or $t_{\max}$
without loss of generality. This simplifies the formulation by reducing the number of nonlinear terms, thereby benefiting the solution process. We present two equivalent formulations:

\paragraph{Circle Packing Formulation}

By fixing $t_{\min} = 1$, the problem becomes a variant of the uniform circle packing problem where points
maintain a pairwise distance of at least 1, with the maximum distance to be minimized. The optimization problem can be formulated as:
\begin{alignat}{2}
    \min\  t_{\max}&
    \label{eq:maxdist-obj} && \\[0.3em]
    \text{s.t.}\quad
     \sum_{k \in \mathcal{D}} (x_{i,k} - x_{j,k})^2 &\ge 1,
    &&\textrm{\quad for all } i,j \in \mathcal{N},\ i<j
    \label{eq:maxdist-min} \\
     \sum_{k \in \mathcal{D}} (x_{i,k} - x_{j,k})^2 &\le t_{\max},
    &&\textrm{\quad for all } i,j \in \mathcal{N},\ i<j
    \label{eq:maxdist-max} \\
     t_{\max} &\ge 1 &&
\end{alignat}
Constraints \eqref{eq:maxdist-min} ensure that all pairwise squared distances are
at least\footnote{For any (locally) optimal solution some of the inequalities in  \eqref{eq:maxdist-min} will be satisfied with equality, otherwise the solution could be trivially improved by scaling.} 1, while constraints \eqref{eq:maxdist-max} bound them by $t_{\max}$. Minimizing $t_{\max}$ is then equivalent to minimizing the ratio $r = \frac{t_{\max}}{1} = t_{\max}$.

\paragraph{Dual Formulation.}
Alternatively, fixing $t_{\max} = 1$ yields a dual perspective where we maximize
the minimum distance subject to a unit bound on the maximum distance:
\begin{alignat}{2}
    \max\  t_{\min} &
    \label{eq:maxdist-obj-simple} && \\[0.3em]
    \text{s.t.}\quad
     \sum_{k \in \mathcal{D}} (x_{i,k} - x_{j,k})^2 &\ge t_{\min},
    &&\textrm{\quad for all } i,j \in \mathcal{N},\ i<j
    \label{eq:maxdist-min-simple} \\
     \sum_{k \in \mathcal{D}} (x_{i,k} - x_{j,k})^2 &\le 1,
    &&\textrm{\quad for all } i,j \in \mathcal{N},\ i<j \\
     0 \le t_{\min} &\le 1 &&
\end{alignat}
Maximizing $t_{\min}$ is then also equivalent to minimizing the ratio $r = \frac{1}{t_{\min}}$. Again, this simple flip renders a nonlinear rational objective into a linear objective, which is preferable for the solution process.

Both formulations yield nonconvex quadratically constrained optimization problems (QCPs) that
can be solved efficiently by modern global optimization solvers. In our experiments the first model performed slightly better, but the difference was not very significant.

The power of the mathematical modeling approach is that we can handle point configurations in any dimension with the same model. Crucially, and this may be one of the reasons why optimization solvers perform great on this problem, the model is essentially unconstrained (at least in its original form): any finite set of points has a largest and smallest pairwise distance, so their ratio is well-defined. There is no extra constraint that the point set has to satisfy (such as inclusion in a bounding region). The constraints get  introduced only through the auxiliary variable for the objective.

\subsection{Computational Results}
\begin{figure}[ht]
\begin{center}
\includegraphics[width=0.24\linewidth]{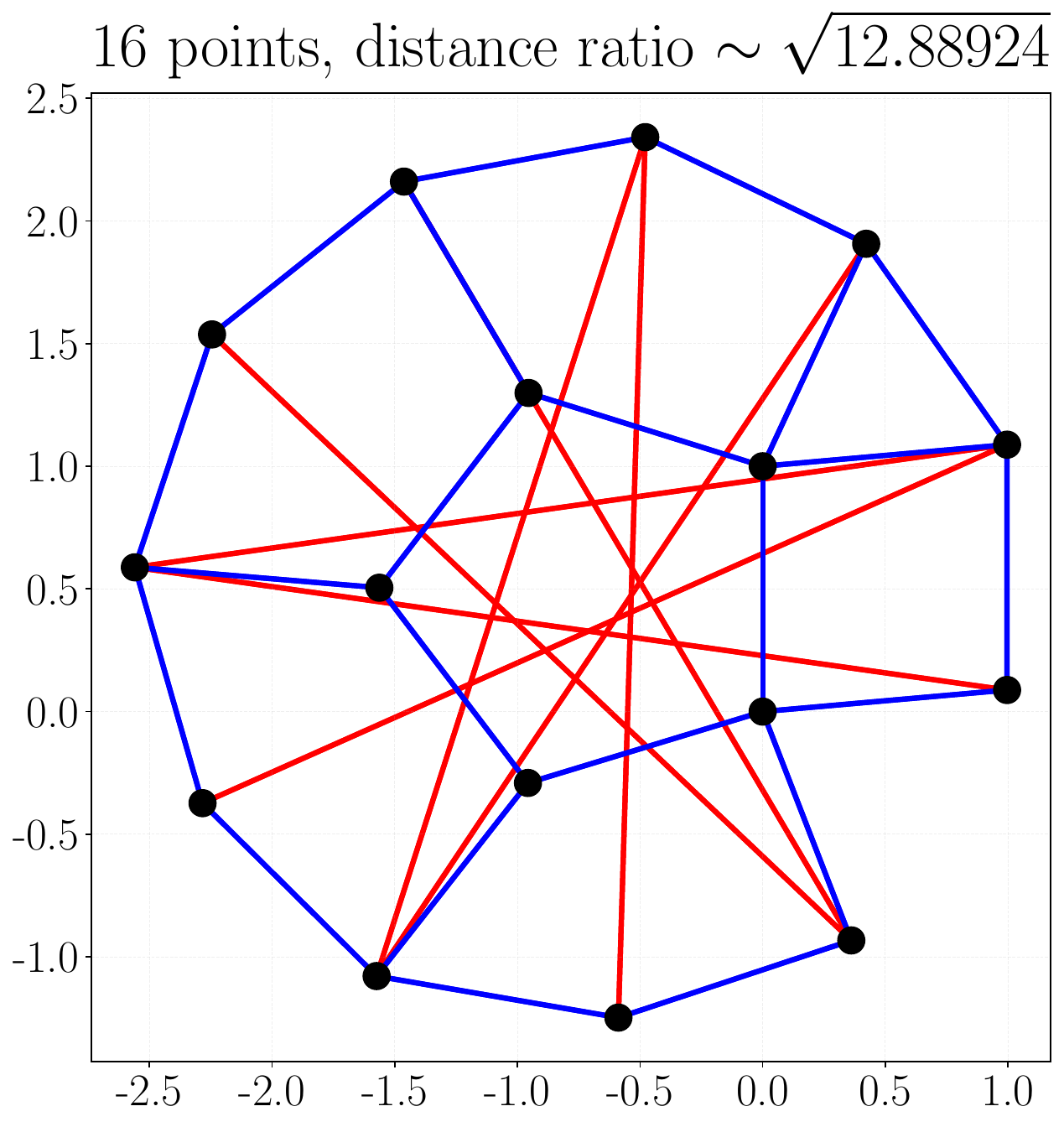}
\includegraphics[width=0.24\linewidth]{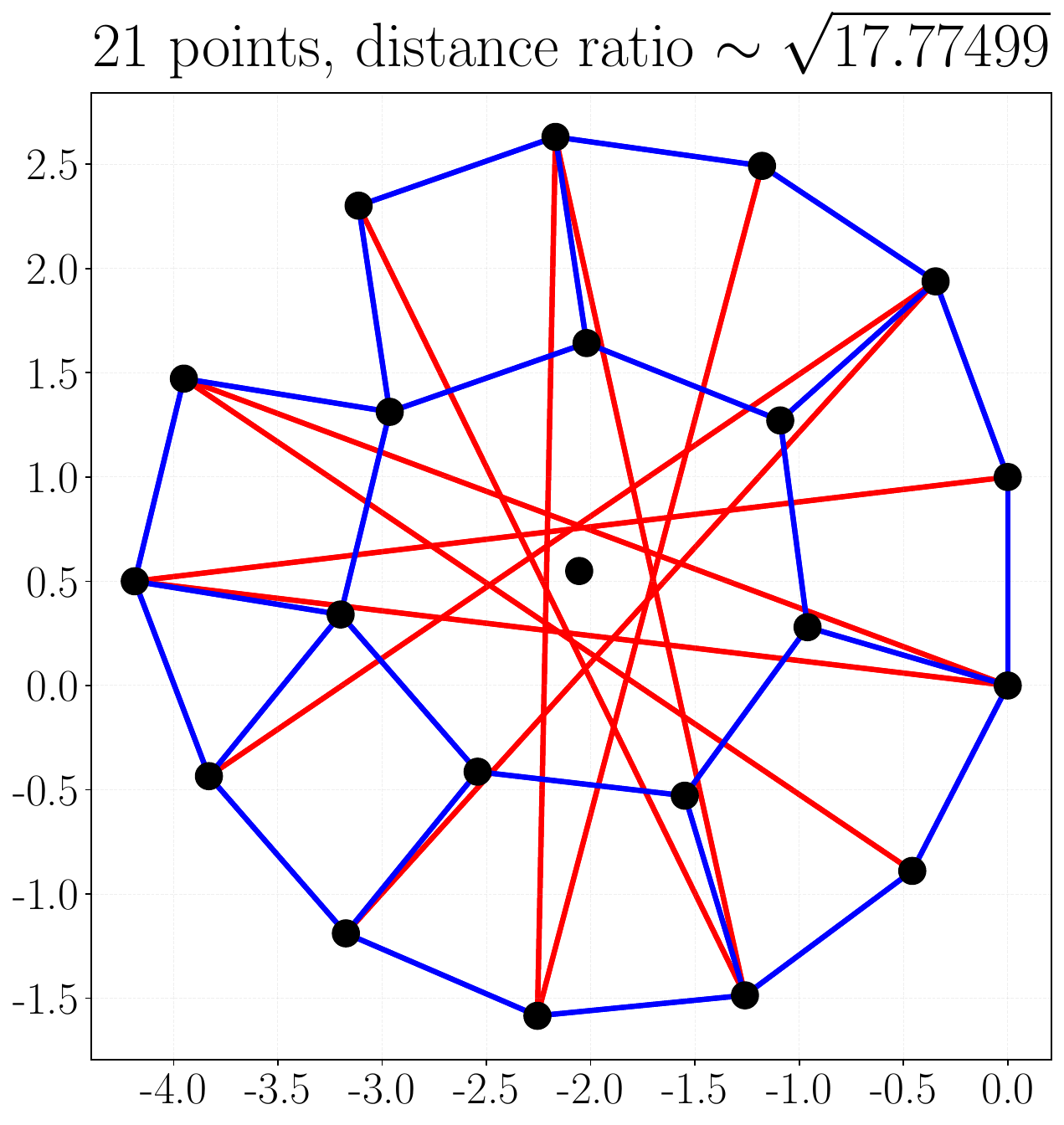}
\includegraphics[width=0.24\linewidth]{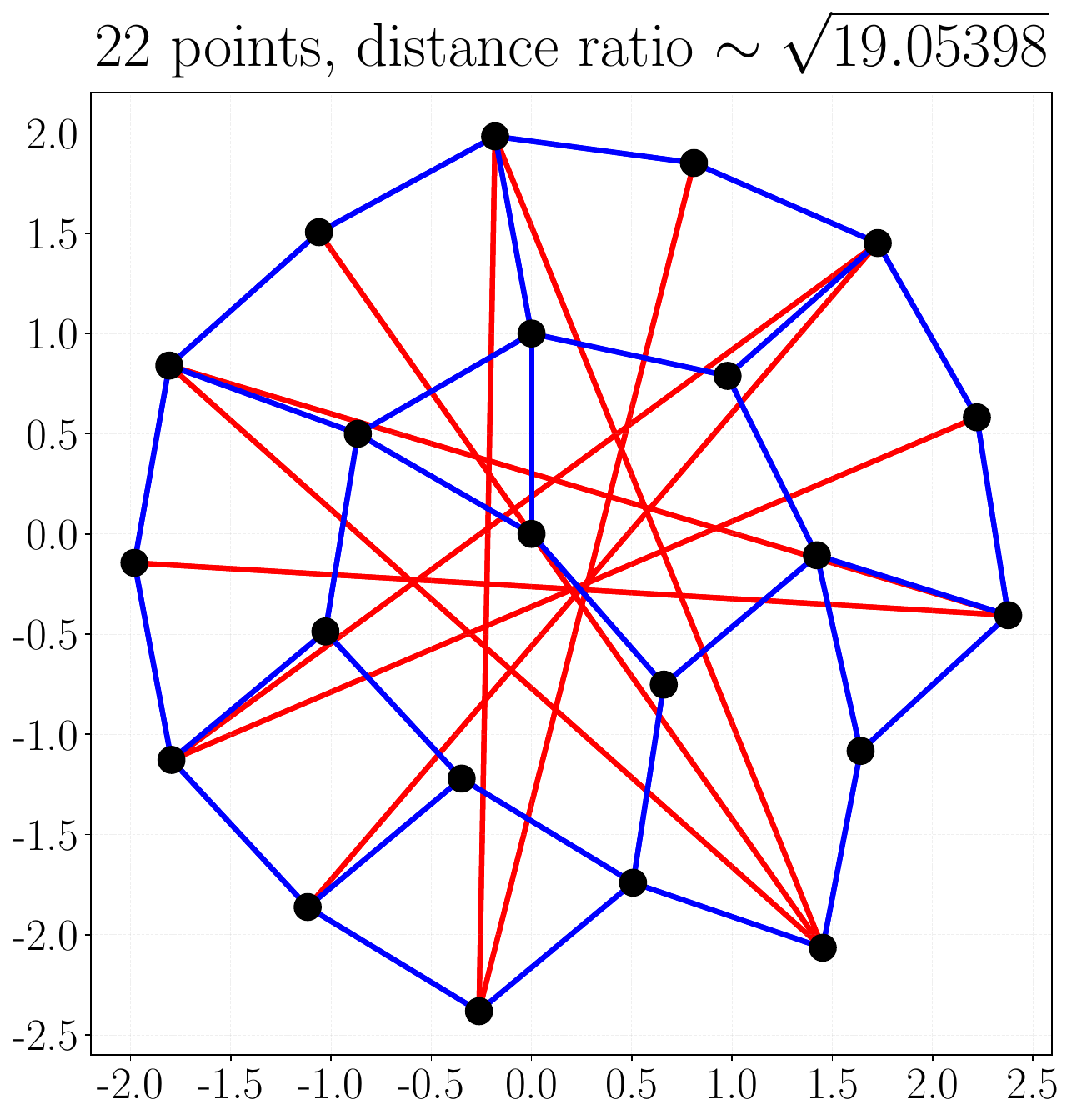}
\includegraphics[width=0.24\linewidth]{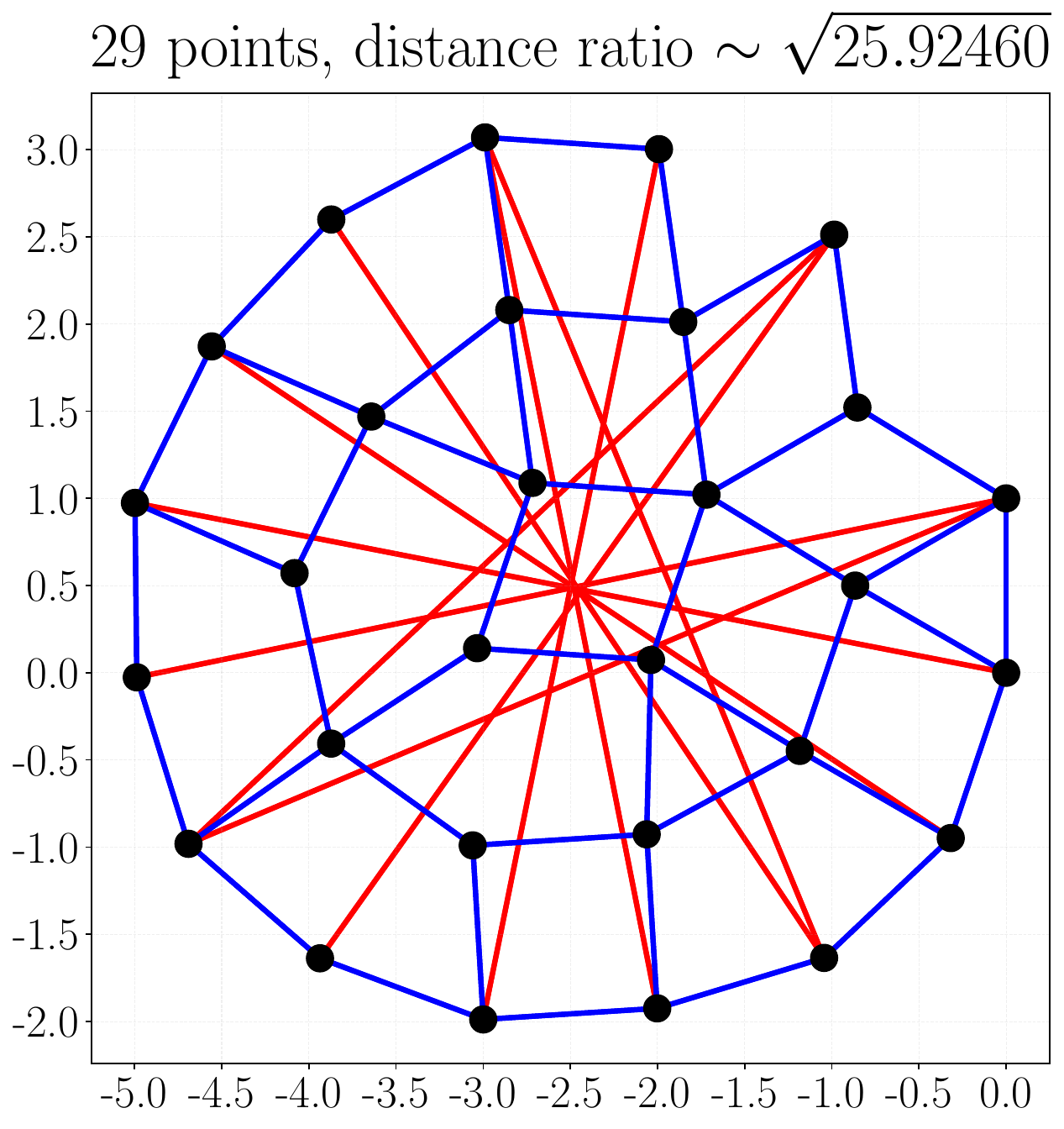}%\\
\end{center}
\caption{Graphical representation of the solutions for the distance ratio problem in 2D. Red lines indicate maximum distance pairs, blue lines indicate minimum distance pairs.}\label{fig:distratio}
\end{figure}

We implemented the above model in the Mosel modeling language \cite{colombani2002mosel}. The resulting instances range in size from 7 variables and 10 constraints (three points in two dimensions) to 91 variables and 876 constraints (30 points in three dimensions).
%All instances were solved using SCIP and FICO Xpress.
The improving solutions (rounded up to five decimal places) reported in Table \ref{tab:distratio} were obtained with this approach, and the corresponding configurations for the 2d case are illustrated in Figure \ref{fig:distratio}. For many additional instances, the method reproduced the currently best-known solutions.
\begin{table}[ht]
\begin{center}
\setlength{\tabcolsep}{8pt}
\begin{tabular}{rc rrr}
\toprule
$d$ &$n$ & Squared distance ratio & Previous best & Source\\
\midrule
2 & 16 & 12.8892\textbf{4} & 12.88927 &\cite{georgiev2025mathematical,novikov2025alphaevolve} \\
2 & 21 & 17.77\textbf{499} & 17.776 & \cite{friedmanPacking} \\
2 & 22 & 19.05\textbf{398} & 19.055 & \cite{friedmanPacking} \\
2 & 29 & 25.92\textbf{460} & 25.929 & \cite{audet2010note} \\
\midrule
3 & 14 &  4.165\textbf{78} &  4.16585 & \cite{georgiev2025mathematical,novikov2025alphaevolve} \\
%2 & 16 & 12.889229767960300 & 12.889266112034630 (AE 2025) \\
%2 & 21 & 17.774980349348699 & 17.776 (David Cantrell 2009) \\
%2 & 22 & 19.053984310699800 & 19.055 (David Cantrell 2009) \\
%2 & 29 & 25.924603726580401 & 25.929 (Charles Audet and Xavier Fournier and Pierre Hansen and Frédéric Messine 2010) \\
%3 & 14 &  4.165783424060720 &  4.165849767225653 (AE 2025) \\
\bottomrule
\end{tabular}
\end{center}
\caption{Improving solutions for the min-max distance ratio problem in 2D and 3D}\label{tab:distratio}
\end{table}

\section{Packing Circles inside a Square or a Rectangle}\label{sec:circlepacking}

The study of packing circles in polygonal regions dates back almost 200 years~\cite{bolyai1832tentamen}.
One of the most studied variants concerns packing a fixed number of unit circles into the smallest possible square; a good survey is given by Peikert \cite{peikert2007packing}. For most instances of up to twenty circles, exact optimal packings are known and can be derived using algebraic geometry techniques based on polynomial systems and Gröbner basis computations, see, e.g.,~\cite{buchberger1976theoretical}.

In this work, we investigate a variant for which it is harder to prove optimality since both the positions and radii of the circles are decision variables. Specifically, we consider packing circles of variable radii into a unit square, or, in a mild relaxation, into a rectangle of perimeter four, while maximizing the sum of their radii.
Best-known solutions are mostly due to Cantrell~\cite{friedman2012circle}, but optimality proofs are only known for some trivial cases.

\subsection{Optimization model}

Given a positive integer $n$, we consider the problem of packing $n$ circles inside a rectangle such that the sum of their radii is maximized. We
examine two variants: packing circles into any rectangle of perimeter 4, and the special case of packing circles into a unit square (which has perimeter 4). A key advantage of mathematical optimization approaches is that the optimization model can be changed easily. This facilitates exploring related problems by simply adding or modifying constraints, such as for the two problem variants shown below.

Let $\mathcal{N} = \{1, 2, \ldots, n\}$ denote the set of circles. We define the
following decision variables:
\begin{description}
\item[$(x_i, y_i) \in \mathbb{R}^2$:] coordinates of the center of circle $i \in \mathcal{N}$
\item[$r_i \in \mathbb{R}_+$:] radius of circle $i \in \mathcal{N}$
\item[$\alpha \in \mathbb{R}_+$:] width of the rectangle
\end{description}
We consider the problem of packing circles in a rectangle with fixed perimeter
$P = 4$ and variable side lengths. We introduce a decision variable $\alpha$
representing the width and set the height to $H = \frac{P}{2} - \alpha = 2 - \alpha$.
We can assume without loss of generality that $\alpha \leq 1$ (i.e., $\alpha$ is
the shorter side). The optimization problem can be formulated as:
\begin{alignat}{2}
    \max\  \sum_{i \in \mathcal{N}} r_i &
    \label{eq:circle-obj} && \\[0.3em]
    \text{s.t.}\quad
     r_i \le x_i &\le \alpha - r_i,
    &&\textrm{\quad for all } i \in \mathcal{N}
    \label{eq:circle-bound-x} \\
    r_i \le y_i &\le (2-\alpha) - r_i,
    &&\textrm{\quad for all } i \in \mathcal{N}
    \label{eq:circle-bound-y} \\
     (x_i - x_j)^2 + (y_i - y_j)^2 &\ge (r_i + r_j)^2,
    &&\textrm{\quad for all } i,j \in \mathcal{N},\ i<j
    \label{eq:circle-nonoverlap} \\
     0 \le r_i &\le \tfrac{\alpha}{2},
    &&\textrm{\quad for all } i \in \mathcal{N}
    \label{eq:circle-radius-bound} \\
    0 < \alpha &\le 1
    \label{eq:circle-width} &&
\end{alignat}
The objective \eqref{eq:circle-obj} is to maximize the sum of all radii.
Note that this is fundamentally different from maximizing the total area covered by circles,
which would be $\pi \sum_{i \in \mathcal{N}} r_i^2$, whereas the linear objective tends
to favor more balanced distributions of circle sizes.
Constraints \eqref{eq:circle-bound-x} and \eqref{eq:circle-bound-y} ensure that
each circle remains entirely within the rectangle: the center coordinates must
maintain a distance of at least $r_i$ from all rectangle boundaries.
Constraints \eqref{eq:circle-nonoverlap} ensure that no two circles overlap by
requiring that the Euclidean distance\footnote{Similar to the distance ratio problem, it is again advantageous to work with squared distances.} between any pair of circle centers is at
least the sum of their radii. Constraints \eqref{eq:circle-radius-bound}
provide an upper bound on each radius: no circle can have a diameter exceeding
the width of the rectangle $\alpha$ (which is the shorter side by assumption).
Finally, constraint \eqref{eq:circle-width} bounds the width variable.

The aspect ratio determined by $\alpha$ is a decision variable that can be modified to maximize
the sum of radii for a given number of circles. We can trivially change this formulation to packing into a unit square by fixing $\alpha = 1$. This is a crucial property of mathematical optimization modeling: the user needs to change only the model and does not have to worry about whether this changes the algorithm: the solvers will take care of that. This contrasts with many heuristic approaches, including LLM-generated ones, in which a new set of heuristics often needs to be developed once the model formulation changes.

Again, we have only linear and nonconvex quadratic constraints.

\subsection{Computational Results}

\begin{figure}[ht]
\begin{center}
\includegraphics[width=0.27\linewidth]{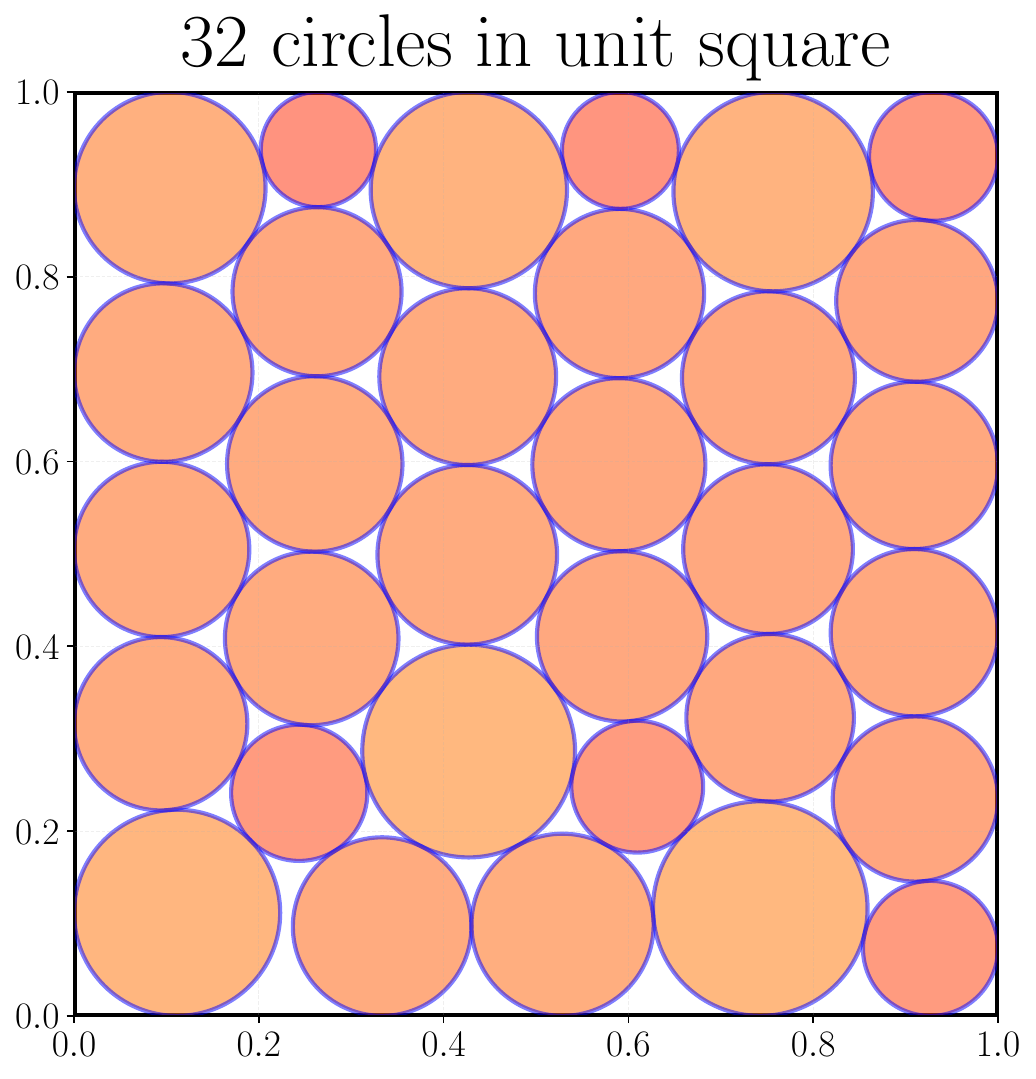}
\includegraphics[width=0.3\linewidth]{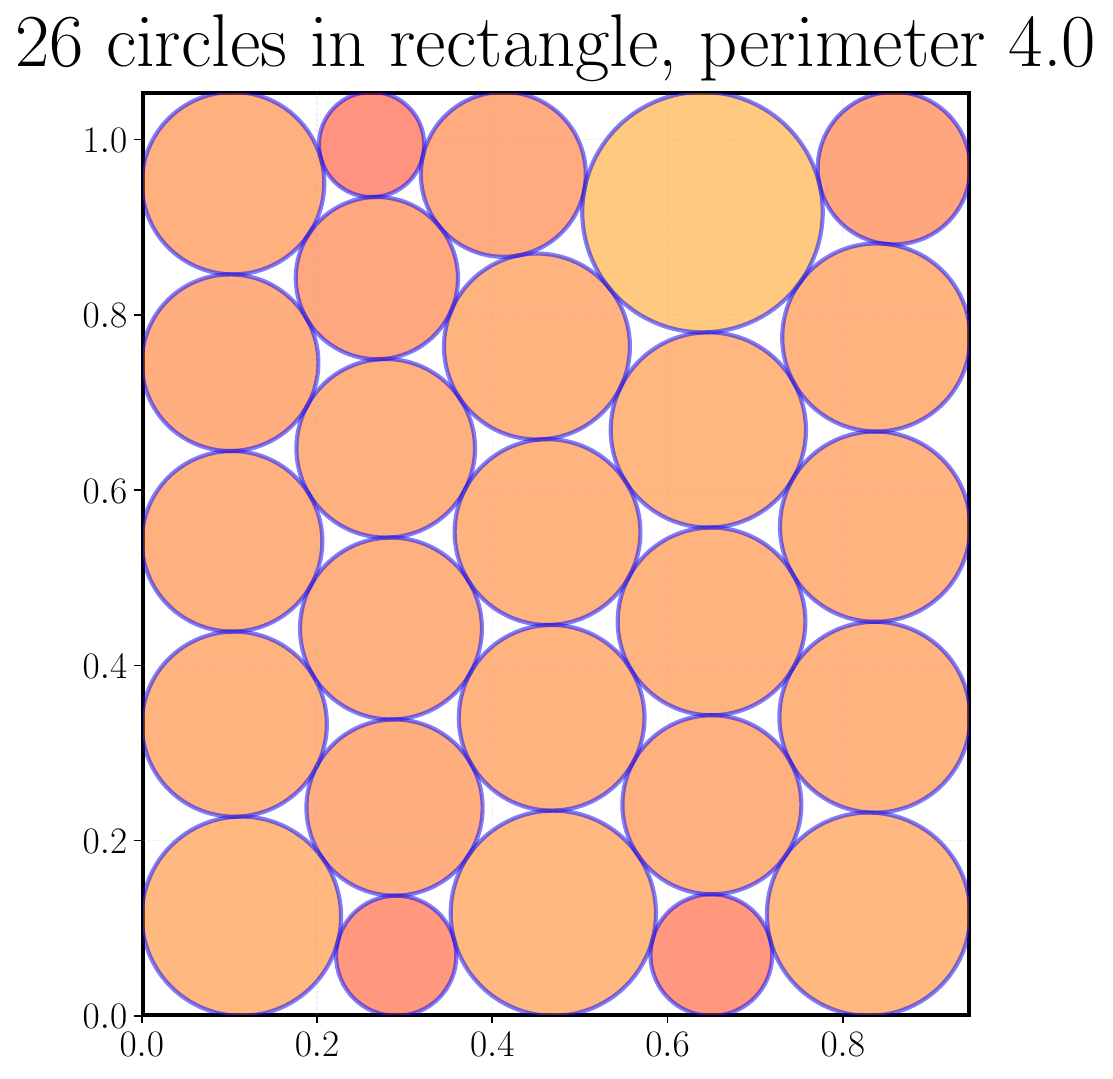}
\includegraphics[width=0.3\linewidth]{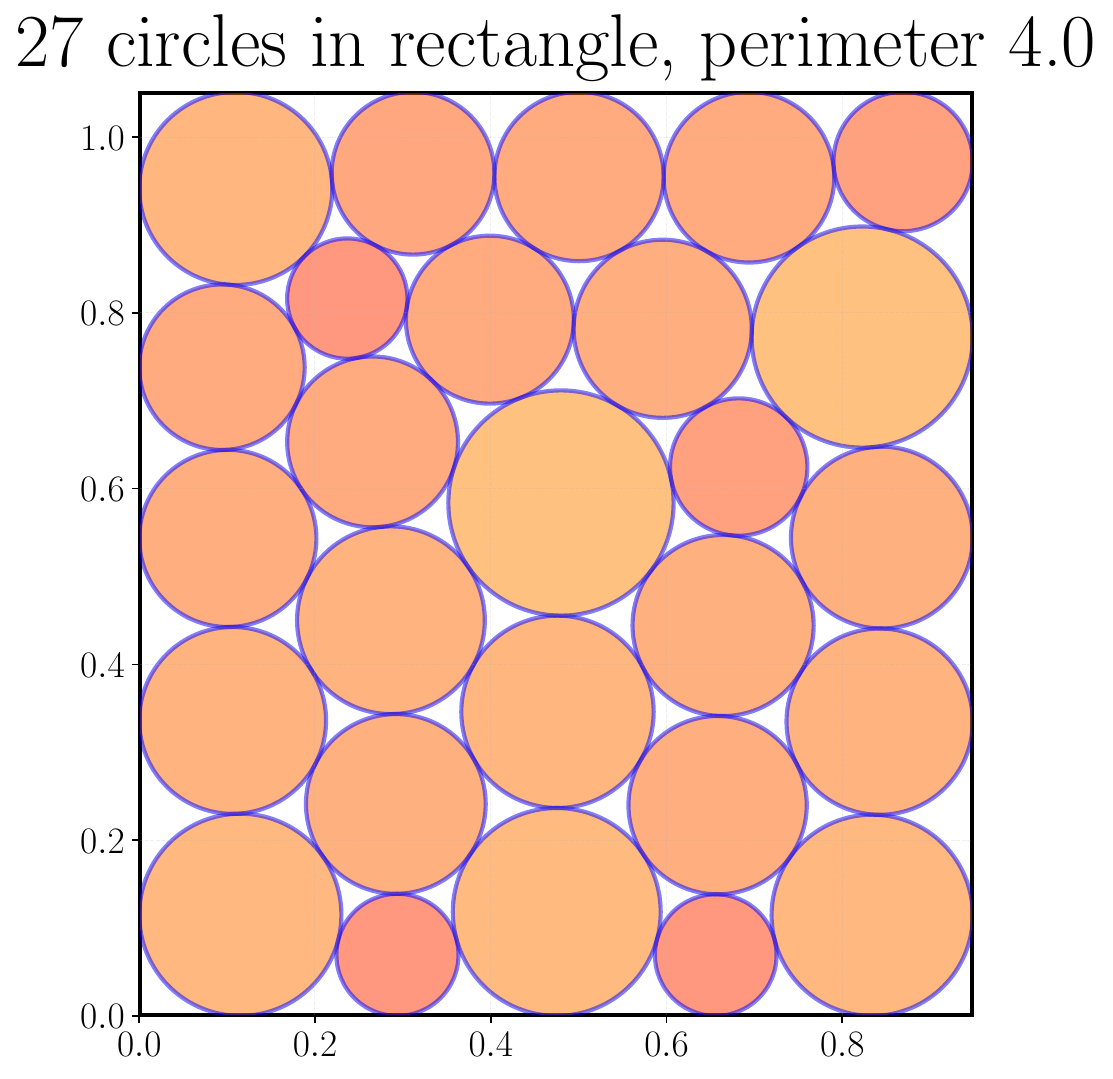}
\end{center}
\caption{Graphical representation of the solutions for the circle packing problem. Left: square variant (n=32). Middle and right: rectangle variant (n=26, n=27).}\label{fig:circlepacking}
\end{figure}

First, let us focus on the restricted problem, packing circles into a square. The following improving solution (floored to five decimal digits) has been found\footnote{We have found possibly even more improving solutions, but the original source \cite{friedman2012circle} for the best-known packings lists only three digits after the decimal. In almost all cases, we managed to match the best-known solutions from literature. } with both solvers (see Table~\ref{tab:packcircle}). The solution is visualized in Figure~\ref{fig:circlepacking}.
\begin{table}[ht]
\begin{center}
\setlength{\tabcolsep}{8pt}
\begin{tabular}{ccrrr}
\toprule
Variant &$n$ & Sum of radii & Previous best & Source\\
\midrule
square &32 & 2.93\textbf{957} & 2.93794 &\cite{georgiev2025mathematical,novikov2025alphaevolve} \\
% 32 & 2.939557549912300 & 2.937944526205518 (AE 2025) \\
\midrule
rectangle & 26 & 2.63\textbf{930} & 2.638 & \cite{friedmanPacking} \\
rectangle & 27 & 2.6\textbf{9015} & 2.687 & \cite{friedmanPacking} \\
% 26 & 2.639308122181169 & 2.638 (David Cantrell 2011) \\
% 27 & 2.690155081571631 & 2.687 (David Cantrell 2011) \\
\bottomrule
\end{tabular}
\end{center}
\caption{Improving solutions for the circle packing problem}\label{tab:packcircle}
\end{table}

Now we can turn to the slightly relaxed version, where instead of a square we are trying to pack into a rectangle of perimeter 4. Obviously, all the previous solutions are still feasible for this relaxed problem, but we can do slightly better.
We found the following improving solutions, again 
%have been found with SCIP~10 by its multistart heuristic within 2 hours
floored to five decimal digits (see Table~\ref{tab:packcircle}). The solutions are visualized in Figure~\ref{fig:circlepacking}. Not surprisingly, the solutions all tend to use a rectangle that is almost square.

\section{Packing Unit Regular Hexagons inside a Regular Hexagon}
\label{hexhex}

Packing various shapes into other shapes in an optimal way has a vast literature, from classical geometry problems and puzzles~\cite{bezdek2010classical,conway2013sphere}, to packing problems solved for designing packaging and packing boxes into shipping containers~\cite{bortfeldt2013constraints,wascher2007improved}.
In general, the more complicated the shape to be packed is, the harder the problems become both to analyze theoretically and to solve practically.
Problems with unit circles are well-studied (compare Section~\ref{sec:circlepacking}), but more complex shapes are often only handled with ad-hoc heuristics.

In this section we consider the problem of packing $n$ regular hexagons with unit side
length into a regular hexagon of minimum side length $R$. Each inner hexagon can be positioned and rotated freely. This
problem is significantly more involved than circle packing due to the possible rotation of the
non-circular geometry of hexagons. However, it is a special case of the polygon packing problem,
where two-dimensional geometrical objects bounded by an arbitrary closed sequence of
non-intersecting line segments must be packed into a container in some efficient way. An
optimization model for the polygon packing problem with a fixed size container was developed
by \cite{kallrath2009cutting}, which was later generalized to convex objects
by \cite{stoyan2016quasi} and afterwards even further to also allow nonconvex objects
by \cite{romanova2018packing}. These approaches construct so-called quasi-phi-functions for certain
object geometries to formulate pairwise non-overlapping conditions among them. The following
hexagon packing formulation adopts this approach in a simplified way by defining Farkas multipliers (see, e.g.,~\cite{bertsimas1997introduction})
on pairwise hexagon intersections in order to exclude interior overlaps, respectively.

\subsection{Optimization model}

Let $\mathcal{N} = \{1, 2, \ldots, n\}$ denote the set of inner hexagons and
$\mathcal{S} = \{0, 1, 2, 3, 4, 5\}$ denote the six sides of any hexagon. We use
the following geometric constants:
\begin{description}
\item[$\rho = \frac{\sqrt{3}}{2}$:] inradius (distance from center to edge, also called apothem)
\item[$\phi = \frac{\pi}{3}$:] angular separation between adjacent vertices
\end{description}
We define the following decision variables:
\begin{description}
\item[$R \in \mathbb{R}_+$:] side length of the outer hexagon (to be minimized)
\item[$(x_i, y_i) \in \mathbb{R}^2$:] coordinates of the center of inner hexagon $i \in \mathcal{N}$
\item[$\theta_i \in \lbrack0, \phi\rbrack$:] rotation angle of inner hexagon $i \in \mathcal{N}$
\end{description}
To facilitate the formulation of the geometric constraints we introduce:
\begin{description}
\item[$a_{i,j}, b_{i,j} \in \mathbb{R}$:] normal vector components for side $j \in \mathcal{S}$ of hexagon $i \in \mathcal{N}$
\item[$c_{i,j} \in \mathbb{R}$:] offset term for side $j \in \mathcal{S}$ of hexagon $i \in \mathcal{N}$
\item[$\lambda_{i,j,k} \in \mathbb{R}_+$:] Farkas multipliers for separation of hexagons $i, j \in \mathcal{N}$ with $i < j$, where $k \in \{1, \ldots, 12\}$
\end{description}
The optimization problem can be formulated as:
\allowdisplaybreaks
\begin{alignat}{2}
    \min\  R& \label{eq:hex-obj} && \\[0.3em]
    \text{s.t.}\quad
     R\rho + \sin(k\phi)\bigl(x_i + \sin(\theta_i + (j + 0.5)\phi)\bigr) \nonumber &\\
     + \cos(k\phi)\bigl(y_i + \cos(\theta_i + (j + 0.5)\phi)\bigr)
      &\ge 0,
    &&\forall i \in \mathcal{N},\ j,k \in \mathcal{S}
    \label{eq:hex-containment} \\[0.3em]
     a_{i,j} = \sin(\theta_i + j\phi),&
    &&\forall i \in \mathcal{N},\ j \in \mathcal{S}
    \label{eq:hex-normal-x} \\
     b_{i,j} = \cos(\theta_i + j\phi),&
    &&\forall i \in \mathcal{N},\ j \in \mathcal{S}
    \label{eq:hex-normal-y} \\
     c_{i,j} = a_{i,j}x_i + b_{i,j}y_i - \rho,&
    &&\forall i \in \mathcal{N},\ j \in \mathcal{S}
    \label{eq:hex-offset} \\[0.3em]
     \sum_{k=1}^{12} \lambda_{i,j,k} &= 1,
    &&\forall i,j \in \mathcal{N},\ i<j
    \label{eq:hex-farkas-sum} \\
     \sum_{k=0}^{5} \lambda_{i,j,k+1} a_{i,k}
      + \sum_{k=0}^{5} \lambda_{i,j,k+7} a_{j,k} &= 0,
    &&\forall i,j \in \mathcal{N},\ i<j
    \label{eq:hex-farkas-x} \\
     \sum_{k=0}^{5} \lambda_{i,j,k+1} b_{i,k}
      + \sum_{k=0}^{5} \lambda_{i,j,k+7} b_{j,k} &= 0,
    &&\forall i,j \in \mathcal{N},\ i<j
    \label{eq:hex-farkas-y} \\
     \sum_{k=0}^{5} \lambda_{i,j,k+1} c_{i,k}
      + \sum_{k=0}^{5} \lambda_{i,j,k+7} c_{j,k} &\ge 0,
    &&\forall i,j \in \mathcal{N},\ i<j
    \label{eq:hex-farkas-sep} \\[0.3em]
     0 \le \theta_i &\le \phi,
    &&\forall i \in \mathcal{N}
    \label{eq:hex-rotation} \\
     \lambda_{i,j,k} &\ge 0,
    &&\forall i,j,k \\
     R &\ge \sqrt{n}
    \label{eq:hex-bounds}
    \end{alignat}

Constraints \eqref{eq:hex-containment} ensure that all vertices of each inner
hexagon lie within the outer hexagon. For a regular hexagon centered at $(x_i, y_i)$
with rotation angle $\theta_i$, the vertices are located at positions
$(x_i + \sin(\theta_i + (j + 0.5)\phi), y_i + \cos(\theta_i + (j + 0.5)\phi))$
for $j \in \mathcal{S}$. Each vertex must satisfy all six half-space constraints
defining the outer hexagon.

Constraints \eqref{eq:hex-normal-x}--\eqref{eq:hex-offset} define the
half-space representation of the inner hexagons. Each side $j \in \mathcal{S}$ of
hexagon $i \in\mathcal{N}$ is characterized by a normal vector $(a_{i,j}, b_{i,j})$ and an
offset $c_{i,j}$, representing the inequality $a_{i,j} x + b_{i,j} y \geq c_{i,j}$.

Constraints \eqref{eq:hex-farkas-sum}--\eqref{eq:hex-farkas-sep} represent
Farkas-based separation conditions to ensure that no two inner hexagons overlap.
These constraints are based on the Farkas lemma for certifying infeasibility of
sets of linear constraints. For each pair of hexagons $(i, j)$ with $i < j$, we introduce
12 Farkas multipliers $\lambda_{i,j,k}$ (one for each side of the two hexagons). The constraints
ensure that the mutual system of inequalities defining the intersection of both hexagons is degenerate, i.e., no point can simultaneously lie in the interior of both
hexagons, thereby guaranteeing non-overlapping configurations.

Constraint \eqref{eq:hex-rotation} restricts the rotation angle to $[0, \phi]$
due to the six-fold rotational symmetry of regular hexagons because rotations beyond
$\phi$ are equivalent to rotations within this range.
Finally, the lower bound on $R$ \eqref{eq:hex-bounds} is implied by elementary geometry. The area of the outer hexagon must be at least the total area of the $n$ contained unit hexagons.

Due to the trigonometric functions, these problems are general nonconvex nonlinear. Even worse, there are nonlinear equality constraints in the formulation we presented. This  can pose a challenge to any solution approach, as equations are naturally harder to satisfy. Global solvers handle this challenge in various ways (both Xpress and SCIP use an outer approximation method).
At the same time, the stronger constrainedness makes this problem also more challenging to create ad-hoc heuristics for.
This is the area where we can report the largest improvements over the previous best solutions, indicating that modern global solvers are well positioned for this strongly constrained problem -- a situation that we are used to from mixed-integer linear optimization for years.

The solution approach can easily be adapted to pack other regular polygons (triangles, squares,
pentagons, etc.) by just changing $n$, recomputing $\rho, \phi$, and rerunning the solver. Actually, any convex polyhedral objects can be handled in a structurally similar way even in higher dimensional domains as long as the linear outer formulations are known because the non-overlapping constraints solely rely on the definition of feasible Farkas multipliers for the intersections. This is one example that demonstrates the generalizable nature of mathematical optimization.

\subsection{Computational Results}

\begin{figure}[ht]
\begin{center}
\includegraphics[width=0.26\linewidth]{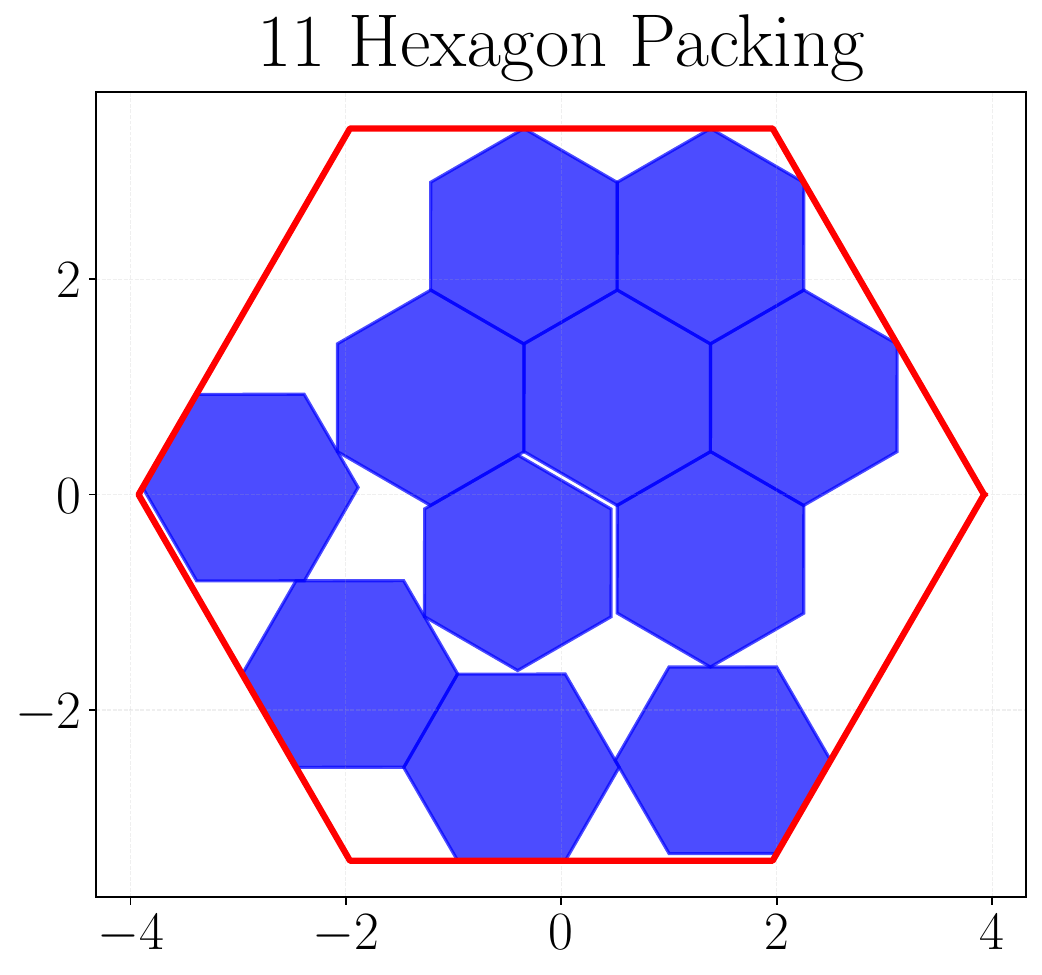}
\includegraphics[width=0.26\linewidth]{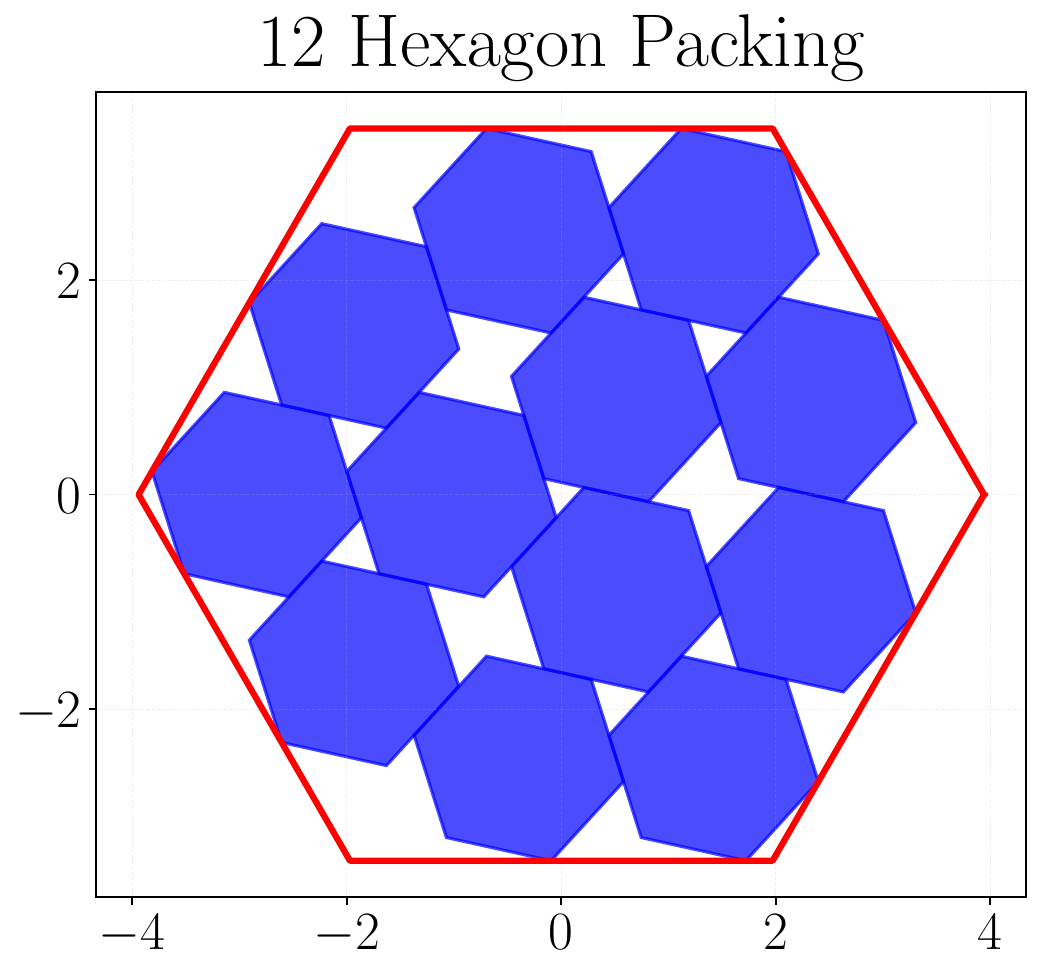}
\includegraphics[width=0.26\linewidth]{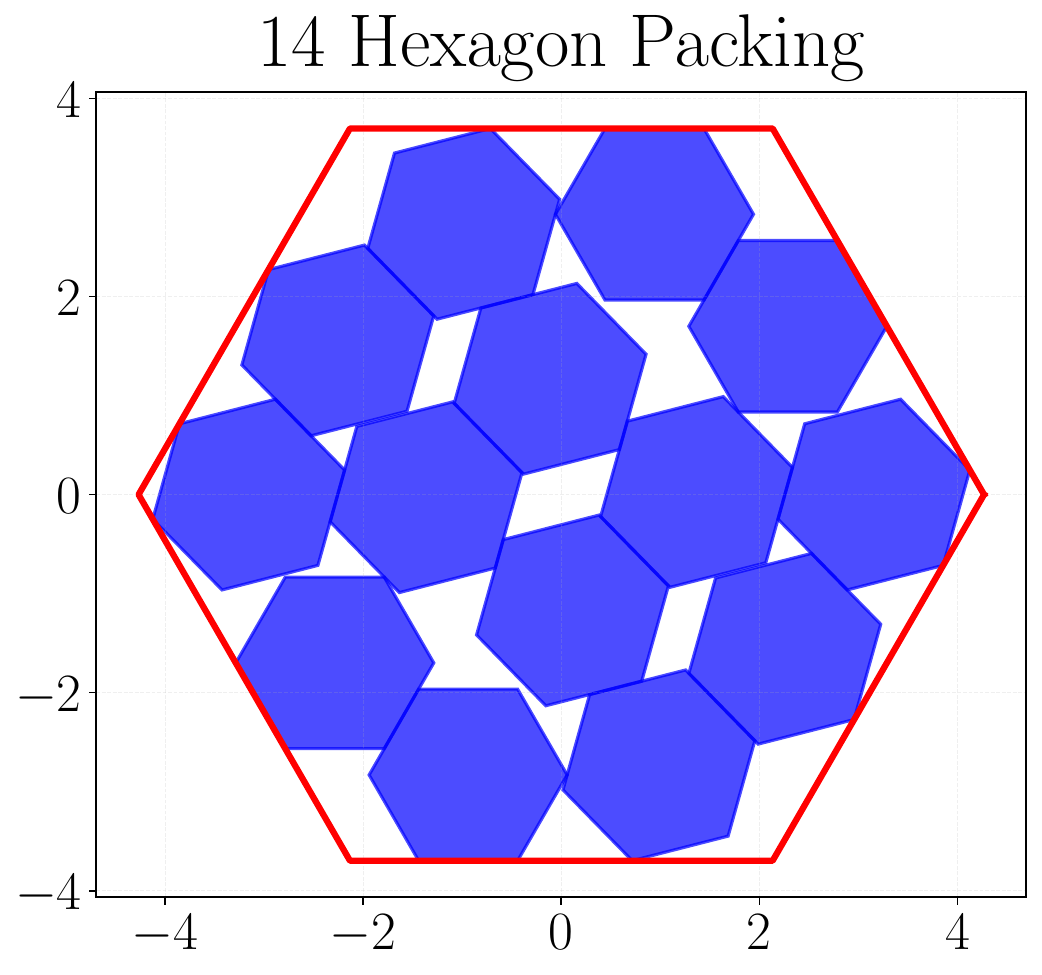}\\
\includegraphics[width=0.26\linewidth]{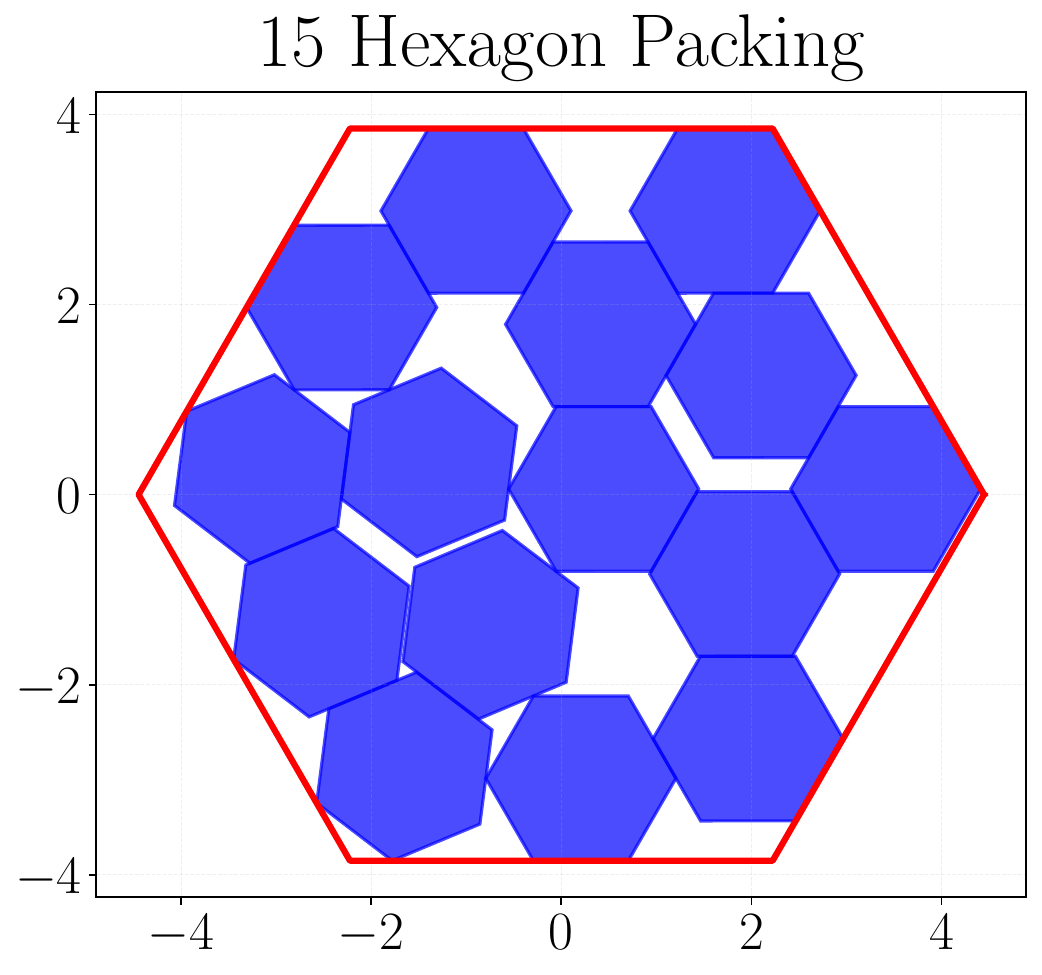}
\includegraphics[width=0.26\linewidth]{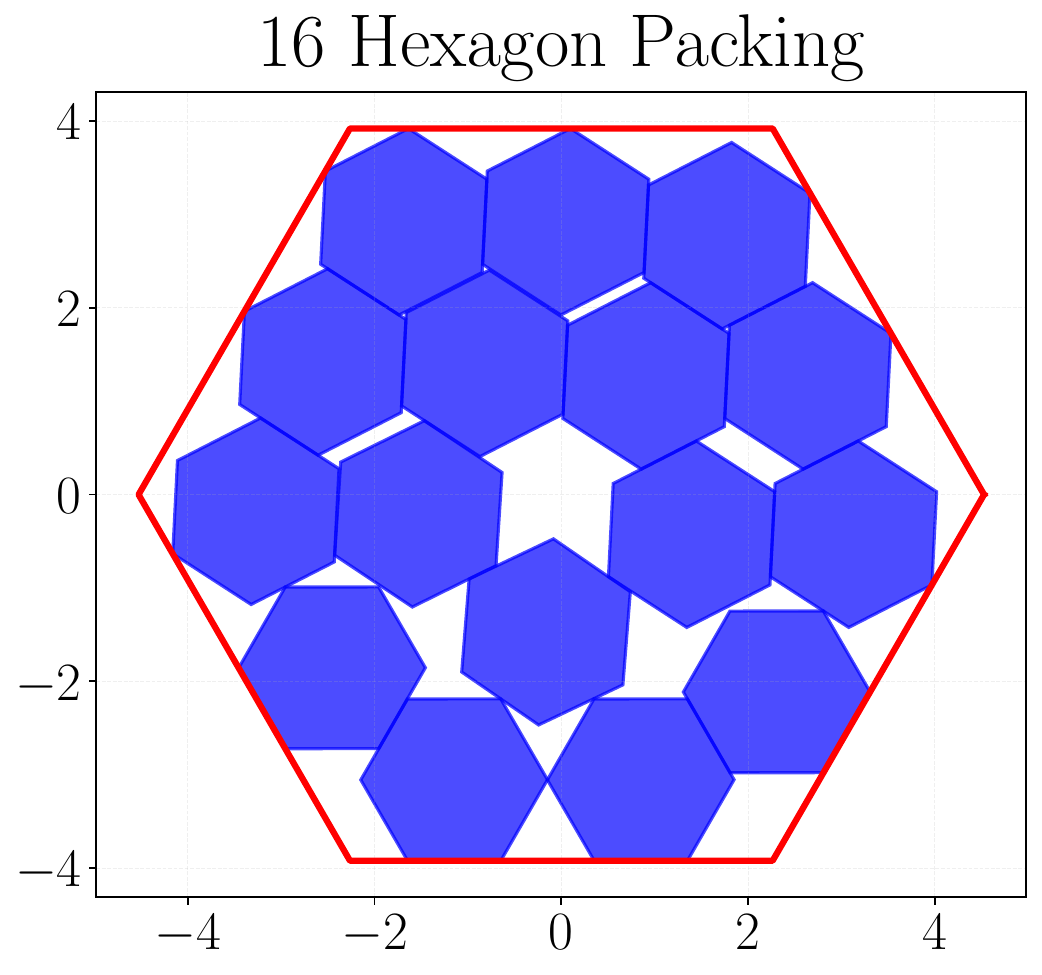}
\end{center}
\caption{Graphical representation of the solutions for the hexagon packing problem}\label{fig:hex}
\end{figure}

%Again, the optimization model was entered into FICO Xpress and SCIP.
The following solutions were produced (often within minutes) with a run using only default settings. The solvers managed to match the best-known solutions for up to 10 hexagons. 
After that, we found solutions that were better than the best-known solutions, except for 13 hexagons, where the known trivial (and most likely, optimal) solution with side length 4 was reproduced.
In particular, the following improving solutions have been found, ceiled to five decimal digits (see Table~\ref{tab:hex} and Figure~\ref{fig:hex}). 
%Similar experiments with SCIP~10 resulted also in better solutions for 11 and 12 hexagons.
%, respectively found within 5 minutes and 2 seconds by its heuristics multistart and subnlp.
\begin{table}[ht]
\begin{center}
\setlength{\tabcolsep}{8pt}
\begin{tabular}{c|rrr}
\toprule
$n$ & Side of containing hexagon & Previous best & Source\\
\midrule
11 & 3.9\textbf{2485} & 3.93010 &\cite{georgiev2025mathematical,novikov2025alphaevolve} \\
12 & 3.941\textbf{65} & 3.94192 &\cite{georgiev2025mathematical,novikov2025alphaevolve} \\
14 & 4.2\textbf{6900} & 4.27240 & \cite{friedman2015hexagon} \\
15 & 4.4\textbf{4769} & 4.45406 & \cite{friedman2015hexagon} \\
16 & 4.5\textbf{2788} & 4.53633 & \cite{friedman2015hexagon} \\
% 11 & 3.924849207777586 & 3.9300911 (AE 2025) \\
% 12 & 3.941642067536949 & 3.9419116 (AE 2025) \\
% 14 & 4.268994957252001 & 4.272391992003231 (Maurizio Morandi 2015) \\
% 15 & 4.447681493753865 & 4.454059275497024 (Maurizio Morandi 2015) \\
% 16 & 4.527878499681161 & 4.536323543632774 (Maurizio Morandi 2015) \\
\bottomrule
\end{tabular}
\end{center}
\caption{Improving solutions for the hexagon packing problem}\label{tab:hex}
\end{table}

\section{Conclusion}
\label{thats_it_folks}
%We have demonstrated that straightforward optimization models coupled with state-of-the-art nonlinear optimization solvers can tackle hard problems in combinatorial geometry. The models are natural and easily modifiable if the problem changes. The runtimes were reasonable (often minutes, typically less than an hour). No solver tuning was performed, both Xpress and SCIP was run with default settings.

Code-generating LLMs provide an effective entry point into new problem domains and ideally facilitate a technology-agnostic exploration of solution approaches.
In the context of the present work, they also helped to draw attention to classical problems where further improvements are still possible and that had not been revisited by mathematical optimization approaches for some time.
Our results show that off-the-shelf global optimization software can solve straightforward formulations of hard combinatorial geometry problems within seconds or minutes.
This highlights the substantial progress achieved in recent years and suggests that nonlinear optimization is approaching the role that (mixed-integer) linear optimization has played for at least two decades: a robust, performant, readily available black-box technology.

At the same time, mathematical optimization can meet several expectations often associated with LLM-based approaches.
Combinatorial and other mathematical  problems can be formulated in a very natural way, and changes in the problem definition typically require only minor modifications to the model and are thereby quickly realized.
Solutions can then be obtained quickly, presumably considerably faster than generating, refining, and testing solver-specific code through dozens of LLM-driven code-generation iterations.

Nevertheless, the two technologies should not be viewed as competing, but rather as complementary.
The solution-method-agnostic prototyping capabilities of an LLM-based code generation are a powerful tool for exploring new problems and formulating new solution approaches. Curiously, in a study with OpenEvolve~\cite{openevolve} on a problem of packing 26 circles specifically, the LLM converged to using an optimization approach~\cite{circlesSLSQP}, namely a  Sequential Least Squares Programming formulation with a local solver implemented in SciPy~\cite{SciPySLSQP,virtanen2020scipy}.
Furthermore, the concept of using LLMs as a support tool for the authoring of optimization models, as, e.g., in~\cite{chen2025optichat,chen2025optimind,huang2025orlm} is a very promising one. 

It is noteworthy that in our experiments the smallest improvements were observed for the relatively unconstrained min-max ratio problem, whereas the largest gains were achieved for the highly constrained hexagon packing problem.
Although based on a limited set of examples, this suggests that optimization-based approaches become increasingly competitive the more constrained the problems are.
Since real-world industrial problems are typically characterized by large and diverse constraint sets, these observations further support the view that global optimization has matured into an industry-ready technology.

\begin{credits}
\subsubsection{\ackname}
We are grateful for the contributions of Liding Xu (ZIB) and the FICO Xpress team,  in particular Bruno Vieira and Susanne Heipcke, for their contributions during an early phase of this project.

Research reported in this paper was partially supported through the Research Campus Modal funded by the German Federal Ministry of Education and Research (fund numbers 05M14ZAM,05M20ZBM) and the Deutsche Forschungsgemeinschaft (DFG) through the DFG Cluster of Excellence MATH+.

\end{credits}
\bibliographystyle{splncs04}
\bibliography{references}

\end{document}